\documentclass{article}

\usepackage{xypic, amssymb, amsmath, amsthm, graphicx, psfrag, setspace,
subfigure}

\newcommand{\chapter}

\makeatletter

\newtheorem{thm}{Theorem}[section]
\newtheorem{cor}[thm]{Corollary}
\newtheorem{lem}[thm]{Lemma}
\newtheorem{prop}[thm]{Proposition}

\theoremstyle{definition}
\newtheorem{defn}[thm]{Definition}
\theoremstyle{remark}
\newtheorem{rem}[thm]{Remark}

\numberwithin{equation}{section}
\newtheorem{eg}[thm]{Example}

\newcommand{\into}{\hookrightarrow}

\newcommand{\id}{\mathbf{1}}
\newcommand{\quotes}[1]{\textquoteleft#1'}

\newcommand{\Hom}{\mbox{Hom}}
\newcommand{\End}{\mbox{End}}

\newcommand{\iso}{\stackrel{\sim}{\longrightarrow}}

\newcommand{\Tr}{\mbox{Tr}}
\newcommand{\C}{\mathbb{C}}

\newcommand{\Z}{\mathbb{Z}}
\newcommand{\cO}{\mathcal{O}}

\newcommand{\Perf}{\mbox{\textit{Perf}}}

\renewcommand{\chapter}{%
    \secdef\Chapter\sChapter}%  typing \chapter{title} invokes
\newcommand{\Chapter}[2][?]{%
    \newpage%
    \refstepcounter{chapter}% update the counter
    \addcontentsline{toc}{chapter}% put something in the toc
        {\protect\numberline{\thechapter} #1}%
    {\raggedright%
        \normalfont\rmfamily\Huge\bfseries%\itshape%
        Chapter \thechapter:\par% write Chapter n
        }
    {\raggedright%
        \normalfont\rmfamily\Huge\bfseries%\itshape%
        #2\par}% write the chapter name
    \chaptermark{#1}% update the text which appears in the header
    \bigskip
    \noindent%
    \!\!}

\newcommand{\sChapter}[1]{
    \newpage%
    {\raggedright%
        \normalfont\rmfamily\Huge\bfseries%\itshape%
        #1\par}% write the chapter name
    \chaptermark{#1}% update the text which appears in the header
    \bigskip
    }

\renewcommand{\subsection}{%
    \secdef\Subsection\sSubsection}%   typing \subsection{title} invokes
\newcommand{\Subsection}[2][?]{%
    \refstepcounter{subsection}%
        \vspace{-0.1\baselineskip}%
    {\flushleft%
        \thesubsection\ \ %
        \normalfont\normalsize\scshape%
        #2%
        \par}%
    \addcontentsline{toc}{subsection}% put something in the toc
        {\protect\numberline{\thesubsection} #1}%
    \nopagebreak
    \vspace{0.5\baselineskip}%
    }

\newcommand{\sSubsection}[1]{%
    \vspace{-0.1\baselineskip}%
    {\flushleft%
        \normalfont\normalsize\sffamily#1%
        \par}%
        \nopagebreak
    \vspace{0.5\baselineskip}%
    }

\renewcommand{\subsubsection}{%
    \secdef\Subsubsection\sSubsubsection}%   typing \subsubsection{title} invokes
\newcommand{\Subsubsection}[2][?]{%
    \refstepcounter{subsubsection}%
    \vspace{-0.1\baselineskip}%
    {\flushleft%
        \thesubsubsection\ \ %
        \normalfont\normalsize\itshape%
        #2%
        \par}%
    \addcontentsline{toc}{subsection}% put something in the toc
        {\indent\protect\numberline{\thesubsubsection} #1}%
    \nopagebreak
    \vspace{0.5\baselineskip}%
    }

\newcommand{\sSubsubsection}[1]{%
    \vspace{-0.1\baselineskip}%
    {\flushleft%
        \normalfont\normalsize\sffamily#1%
        \par}%
        \nopagebreak
    \vspace{0.5\baselineskip}%
    }

\renewcommand{\section}{\@startsection
    {section}% the name.
    {1}% the level.
    {0mm}% the indent before the heading itself (minus means no indent).
    {-0.6\baselineskip}% the space before the heading,
    {0.5\baselineskip}% and the space after (negative is running header).
    {\normalfont\large\scshape\centering}}% the style.

\makeatother

\begin{document}

\title{The closed state space of affine Landau-Ginzburg B-models}%
\author{Ed Segal}%

\maketitle

\begin{abstract}
We study the category of perfect cdg-modules over a curved algebra, and in particular the category of B-branes in an affine Landau-Ginzburg model. We construct an explicit chain map from the Hochschild complex of the category to the closed state space of the model, and prove that this is a quasi-isomorphism from the Borel-Moore Hochschild complex. Using the lowest-order term of our map we derive Kapustin and Li's formula for the correlator of an open-string state over a disc.

\end{abstract}

\tableofcontents
% ----------------------------------------------------------------

\section{Introduction}\label{sectintro}

A Landau-Ginzburg model is a 2-dimensional supersymmetric quantum field theory depending on a K\"ahler manifold $X$ and a holomorphic function $W$ on $X$, called the superpotential. If $X$ is Calabi-Yau then the theory admits a B-twist, which makes it into a topological theory that depends only on the complex structure of $X$ and not the metric. 

In the simplest kind of 2d topological field theory the worldsheet is just a topological 2-manifold, and since there are not very many of these a theory like this cannot contain much information. A Landau-Ginzburg B-model is a more complicated kind of theory where the worldsheet is a Riemann surface, but we integrate over families of complex structures, so we pick up the topology of the moduli space of Riemann surfaces. A theory like this is called mathematically a Cohomological Field Theory or a Topological Conformal Field Theory. 

Thus the physics predicts that given a complex manifold $X$ and a holomorphic function $W$ we should be able to construct a TCFT. This problem is made easier by the results of Costello \cite{costello2} on the formal structure of TCFTs. He showed that the open sector is specified by a choice of Calabi-Yau $A_\infty$-category (or dg-category), called the category of branes. Furthermore if we know the open sector then there is a canonical choice of closed sector, it's the Hochschild complex of the category of branes. 

When $W=0$ it is well known that the category of B-branes is the described by the derived category $D^b(X)$ of coherent sheaves on $X$, or more accurately by a dg- or $A_\infty$- enhancement of it, such as the category $\Perf(X)$ of perfect complexes on $X$. Also, the Hochschild homology of $\Perf(X)$ is the Dolbeault cohomology of $X$, which is the physically-predicted closed sector.

When $W\neq 0$, it was suggested by Kontsevich that B-branes are \quotes{twisted complexes} of vector bundles, i.e. instead of carrying differentials they carry endomorphisms $d$ such that $d^2 = W$. These are related to classical objects from algebraic singularity theory called matrix factorizations, but they were first studied mathematically as B-branes by Orlov \cite{orlov}. In Orlov's work B-branes form a $\Z_2$-graded triangulated category, we will instead be using the construction from $\cite{segal}$ which gives a $\Z$-graded dg-category $Br(X,W)$. 

When $X$ is affine and $W$ has isolated singularities, the physics predicts that the closed sector of the B-model is the Jacobi ring $J_W$. It is hence natural to conjecture that this is the Hochschild homology of $Br(X,W)$ in this case. This result has proved by Dyckerhoff \cite{dyckerhoff}, by identifying a generating object for the category. 

In this paper, we consider the case that $X$ is affine, but $W$ is arbitrary. The natural replacement for $J_W$ is the chain complex
$$(\Omega^\bullet_X, \wedge dW)$$
which we call the \quotes{off-shell closed state space}, when $W$ has isolated singularities this has homology $J_W$. When the singularities of $W$ are not isolated one does not expect $Br(X,W)$ to be Calabi-Yau, so the full TCFT structure does not exist. Nevertheless, the Hochschild homology of $Br(X,W)$ should still be equal to the homology of this complex. 

What we achieve in this paper is the construction of an explicit chain map (\ref{higherbb}) between the Hochschild complex of $Br(X,W)$ and this off-shell closed state space. Unfortunately we cannot prove directly that our map is a quasi-isomorphism. However, our map is naturally defined on a completed\footnote{More accurately \quotes{uncocompleted}.} version of the Hochschild complex called the Borel-Moore Hochschild complex \cite{CT}, and we show our map \textit{is} a quasi-isomorphism on this slightly larger complex. Furthermore, it should follow from the work of Polischuk and Positselski \cite{PolPos} that the two kinds of Hochschild complex are in fact quasi-isomorphic.

When $W$ does have isolated singularities there is a residue map on $\Omega^\bullet_X$ that sends forms to their residues at the singularities. If we apply our chain map to a single morphism in $Br(X,W)$, and then take the residue, we recover the Kapustin-Li formula for the correlator of an open-string state over a disc \cite{KL}. 

It would be nice to generalize to when the space $X$ is not affine. In that case the Hochschild homology of $Br(X,W)$ is presumably the homology of 
$$(\mathcal{A}^{\bullet, \bullet}_X, \;\bar{\partial} + \wedge dW) $$
This was argued physically in \cite{HerbLaz}, and a proof for the $\Z_2$-graded case is sketched in \cite{LinPom}. We do not know how to write down the analogue of our map in the non-affine case, since one must deal with the fact that vector bundles on $X$ can have non-trivial Chern characters.

The outline of this paper is as follows:

In Section 2.1 we discuss \textit{curved algebras}, of which affine Landau-Ginzburg B-models are a special case. We then define the category $\Perf(A,W)$ of perfect curved dg-modules over a curved algebra, which is the more general analogue to the category of B-branes. 

In Section 2.2 we spend some time recalling the non-commutative-geometric language of Kontsevich and Soibelman \cite{KS}, and define the Hochschild and Borel-Moore Hochschild complexes in this language. 

Section 2.3 contains the key idea of this paper. We construct an isomorphism between the Borel-Moore Hochschild complexes of $\Perf(A,W)$ and of another closely related category. This second category has the same objects and morphisms, but it has no differential, instead it has a curvature term induced from $W$. This isomorphism is easy to see in the geometric language, as it comes from a translation of a non-commutative vector space. We can then show that the Borel-Moore Hochschild complex of this curved category is equal to that of the curved algebra.

In Section 3 we discuss Landau-Ginzburg models and their closed state spaces. We apply the results of Section 2, and get some explicit formulas, in particular we recover the Kapustin-Li formula. Finally, we show how to extend these results to affine orbifolds.

\textbf{Acknowledgements.} The first version of this paper erroneously claimed that I had calculated the Hochschild homology of the category of B-branes. I'm grateful to Tobias Dyckerhoff for pointing out the mistake, and to Andrei C\u ald\u araru and Junwu Tu for working out what the correct statement was.

 I'd also like to thank Nils Carqueville, Tom Coates, Kevin McGerty, Ezra Getzler, Kevin Lin, Daniel Pomerleano, Leonid Positselski, and Richard Thomas for many helpful discussions and ideas.

\section{Curved algebras}\label{curvedalgebras}

In this section we prove that the Borel-Moore
Hochschild complex of the category of perfect dg-modules over a curved algebra
is quasi-isomorphic to the Borel-Moore Hochschild complex of the curved algebra itself.

We work over an arbitrary ground field $k$ of characteristic zero, and \textit{category}
means a $k$-linear category.

\subsection{Curved algebras and cdg-modules}\label{cdgasandcdgms}

\begin{defn} \cite{GJ} A \textit{curved dg-algebra} is a triple $(A, d, W)$ where $A$ is a
graded associative algebra, $d$ is a degree 1 derivation of $A$, and $W\in A$ is a degree 2 element such that $dW=0$ and
$$ d^2 = [W, -]$$
A \textit{curved algebra} is a curved dg-algebra where $d=0$, and hence $W$
is central. The definitions of \textit{curved dg-category} and \textit{curved
category} are similar.
\end{defn}
We may choose to work with either a $\mathbb{Z}$-grading or just a $\mathbb{Z}_2$-grading. Of course, in the $\mathbb{Z}_2$-graded version \quotes{$W$ is a degree 2 element} means $W$ is a degree 0 element.

The name comes from thinking of $W$ as the \quotes{curvature} of the \quotes{connection} given by $d$. In this paper we will be mostly concerned with curved algebras
and (non-curved) dg-categories. We will generally drop the $d$, but not the
$W$, from the notation.

\begin{defn} A \textit{cdg-module} over a curved dg-algebra
$(A, W)$ is a pair $(M, d_M)$ where $M$ is a graded $A$-module and $d_M$ is a degree 1 linear endomorphism of $M$ such that
such that
$$d_M(am) = (da)m + (-1)^{|a|}a(d_M m) $$
and
 $$d_M^2=W$$

\end{defn}

If $W=0$ this is just the usual definition of a dg-module over a dga, but having a non-zero $W$ \quotes{twists} the differentials.

Given two cdg-modules $(M, d_M)$ and $(N, d_N)$, we have a graded vector space $ \Hom_A(M,N)$ with a degree 1 endmorphism
\begin{equation}\label{dhom}
d_{M,N}(f):=(d_N\circ f) - (-1)^{|f|}(f\circ d_M)
\end{equation}
As in the case $W=0$, this is in a fact a differential, even though neither $d_M$ nor $d_N$ is. The two copies of $W$ that occur in the square of this expression have opposite signs and cancel. This means that the category of cdg-modules is a dg-category.

\begin{rem}\label{curvedcat}When $A$ is commutative, this definition can be seen as part of a larger structure. Given a cdg-module $M$ over $(A, W_M)$, and another $N$ over $(A, W_N)$, we
can form their tensor product $M\otimes_A N$, which is a cdg-module over $(A, W_M+W_N)$.  This defines a monoidal product on the the category of cdg-$A$-modules
where we allow all possible $W$s, and the $\Hom$ complex defined above is
an internal Hom functor. This larger category is a curved dg-category with
curvature $W_M$ at each object $M$. Arguably even when
considering the sub-category of cdg-modules over a fixed $W$ one should keep
this (now central) curvature term, but we shall not do so in this paper.
\end{rem}

In the ordinary case of a (non-curved) dg-algebra $A$, one often studies not the category of dg-modules but rather the derived category $D(A)$, which is its localization at quasi-equivalences. This is a triangulated category, and we wish to work instead with dg-categories. Now the category of dg-modules over $A$ is of course a dg-category, but it is not the correct dg-category to study since its homotopy category is not $D(A)$. We can correct it by taking only the subcategory of projective modules, this is a dg-category that does have $D(A)$ as its homotopy category. Or we can take some smaller dg-category such as the category of perfect (i.e. finitely-generated projective) dg-modules, which has a slightly different homotopy category that may be better-behaved.

We will work with an analogue of the category of perfect complexes over a  curved algebra.  The derived category of modules over a curved algebra (in fact over a curved $A_\infty$-algebra) has been constructed in \cite{nicolas}, but by analogy with the non-curved case one should not expect it to agree with the homotopy category of our category, except maybe in special cases such as $A$ commutative and smooth.

\begin{defn} A cdg-module $(M, d_M)$ over a curved dga $(A, W)$ is \textit{perfect} if $M$ is a finitely-generated and projective module over the underlying algebra $A$.
\end{defn}

\begin{defn} The dg-category $\Perf(A,W)$ is the full sub-category of the category of cdg-$(A,W)$-modules with objects the perfect cdg-modules.
\end{defn}

\subsection{Some non-commutative geometry}

In this section, which is almost entirely lifted from \cite{KS}, we set up some geometric language for studying curved dg-categories.

\subsubsection{Curved $A_\infty$-structures and polynomial curved $A_\infty$-structures}

Let $V$ be a graded vector space. Initially we'll assume that $V$ is degree-wise finite-dimensional, this is so we can take (graded) duals without worrying. We want to consider $V$ as a \textit{non-commutative vector space}, this means we declare the \quotes{ring of functions} on $V$ to be the ring of non-commutative polynomials
$$\mathcal{O}(V) := TV^\vee = \bigoplus_{k\geq 0} (V^\vee)^{ \otimes k}$$
Many of the usual constructions of algebraic geometry go through unchanged. For example, $\mathcal{O}(V)$ contains an ideal generated by $V^\vee$, and the quotient by the $k$th power of this ideal gives an algebra that corresponds to a $k$th-order neighbourhood of the origin in $V$. If we take the limit of these algebras over $k$ we get the completed tensor algebra
$$\hat{\mathcal{O}}(V) := \hat{T}V^\vee = \prod_{k\geq 0} (V^\vee)^{ \otimes k}$$
which corresponds to a formal neighbourhood of the origin.

A \textit{vector field} on $V$ is a derivation of $\mathcal{O}(V)$, so a derivation of $\hat{\cO}(V)$ is like the germ of a vector field at the origin.

\begin{defn} Let $V$ be degree-wise finite-dimensional. A (curved) $A_\infty$-structure on $V$ is a derivation
$$Q :\hat{\cO}(V) \to \hat{\cO}(V)$$
of degree 1, such that
$$[Q,Q] =0$$
\end{defn}
Since $Q$ is a derivation it is determined by its effect on linear functions, i.e. by the map
$$Q : V^\vee \to \hat{T}V^\vee$$
The components $Q_0, Q_1, Q_2,...$ of this map are the Taylor coefficients of the vector field. Dualizing them gives us an infinite sequence of maps
$$Q^\vee_k : V^{\otimes k} \to V$$
all of degree 1. Now let  let $A=V[-1]$ be the degree-shift of $V$. We get an induced sequence of multilinear maps on $A$, which are conventionally denoted
$$m_k : A^{\otimes k} \to A$$
The map $m_k$ has degree $2-k$. These maps are called the $A_\infty$\textit{-products} on $A$, and $A$ is called an $A_\infty$\textit{-algebra}.

The requirement that $[Q,Q] = 2Q^2 = 0$ translates to relations on the products $m_k$. For example, suppose that the Taylor coefficients of $Q$ vanish above the quadratic term, so $Q=Q_0 + Q_1 + Q_2$. Then $Q^2=0$ iff the following four equations hold:
\begin{eqnarray*} Q_0Q_1 =0\\
Q_1^2 + (Q_0\otimes\id + \id\otimes Q_0)Q_2 = 0 \\
Q_2Q_1 + (Q_1\otimes\id + \id\otimes Q_1)Q_2 =0\\
(Q_2\otimes\id + \id\otimes Q_2)Q_2 =0
\end{eqnarray*}
The $A_\infty$-products on $A$ consist of a degree 0 bilinear product $m_2$, a degree 1 linear endomorphism $m_1$, and  and a degree 2 constant $m_0\in A$. Dualizing the above relations, and inserting the signs required by the degree shift, gives
\begin{eqnarray*} m_1m_0 &=&0\\
m_1^2 &=& m_2(m_0\otimes\id - \id\otimes m_0) \\
m_1m_2&=&m_2(m_1\otimes\id + \id\otimes m_1)\\
m_2(m_2\otimes\id)&=& m_2(\id\otimes m_2)
\end{eqnarray*}
These are precisely the axioms that make $A$ a curved dg-algebra. The curvature is $m_0$, so if we want an ordinary dga then we have to require $Q_0=0$. In fact, in this paper all of our curved $A_\infty$-structures will just be curved dgas, i.e. their Taylor expansions will vanish above the quadratic terms, nevertheless if we use this geometric language then it is simpler to work in this extra generality.

A morphism between two curved $A_\infty$-structures $(V_1, Q_1)$ and $(V_2, Q_2)$ is a homomorphism between the completed algebras
$$\hat{\cO}(V_2) \to \hat{\cO}(V_1)$$
that intertwines the two derivations $Q_2$ and $Q_1$. It follows that an isomorphism between curved $A_\infty$-structures corresponds to pulling-back the germ of a vector field via the germ of a diffeomorphism.

If we have a curved $A_\infty$-structure $(V, Q)$ where the derivation $Q$ has only a finite number of terms, so there are only finitely many non-zero $A_\infty$-products, then $Q$ in fact defines a derivation of $\cO(V)$. This is a vector field defined over all of $V$, not just in a formal neighbourhood of the origin. We shall call such a structure a \textit{polynomial} curved $A_\infty$-structure.

 Similarly a \textit{polynomial morphism} between polynomial curved $A_\infty$-structures $(V_1, Q_1)$ and $(V_2, Q_2)$ is a homomorphism
$$\cO(V_2) \to \cO(V_1)$$
intertwining $Q_2$ and $Q_1$. Polynomial morphisms are not a subset of morphisms, because they need not preserve the origin.

\subsubsection{Differential forms}

Let $V$ be a degree-wise finite-dimensional graded vector space. The \textit{odd tangent bundle} to $V$ is the graded vector space
$$V\oplus V[1]$$
The ring of functions on this is called the space of \textit{de Rham differential forms} on $V$, and denoted
$$ \Omega^\bullet(V) = \bigoplus_{m\geq 0} \Omega^m(V)= \cO(V\oplus V[1])$$
The splitting is by the number of copies of $V[1]^\vee$ that appear, so for example
$$\Omega^0(V) = \cO(V)$$
and
$$\Omega^1(V) = \cO(V) \otimes V[1]^\vee \otimes \cO(V)$$
 Let $\delta_{dR}$ be the linear endomorphism of $V\oplus V[1]$ that maps $V[1]$ isomorphically onto $V$. Then $\delta_{dR}$ has degree $-1$, and squares to zero. It defines a degree 1 linear vector field $d_{dR}$ on $V\oplus V[1]$ which has $\delta_{dR}^\vee$ as its single non-zero Taylor coefficient. Note that
$$d_{dR}: \Omega^m(V) \to \Omega^{m+1}(V)$$
 and $d_{dR}^2=0$.

Let $X$ be any vector field on $V$. We define a vector field $i_X$ on $V\oplus V[1]$ by declaring it to have Taylor coefficients
$$i_X: V[1]^\vee \iso V^\vee \stackrel{X}{\longrightarrow} \cO(V) \into \Omega^\bullet(V)$$
This just contracts a differential form by the vector field $X$. We also define the \textit{Lie derivative} along $X$ as
$$Lie_X = [d, i_X]$$
$Lie_X$  preserves each $\Omega^m(V)$. It is easy to check that
$$[Lie_X, Lie_Y] = Lie_{[X,Y]}$$
 Now let $Q$ be a polynomial curved $A_\infty$-structure on $V$, so $Q$ is degree 1 vector field on $V$ and $[Q,Q]=0$. It follows that $Lie_Q$ is a degree 1 vector field on $V\oplus V[1]$ such that
$$[Lie_Q, Lie_Q] = 2Lie_Q^2 = 0$$

For any graded vector space $V$, the ring of non-commutative polynomials on $V$ carries an action of the infinite cyclic group, by cyclically permuting the factors of each $V^{\vee\otimes k}$. We call the coinvariants the \textit{cyclic functions} on $V$, and denote them by
$$\cO_{cycl}(V)  = (\cO(V))_{\Z} = \cO(V) / [\cO(V), \cO(V)]$$
Any vector field on $V$ preserves $[\cO(V), \cO(V)]$, and hence induces a linear map on $\cO_{cycl}(V)$.

The \textit{cyclic differential forms} on $V$ are
$$\Omega^\bullet_{cycl}(V) := \cO_{cycl}(V\oplus V[1])$$
They carry a differential $d_{cycl}$ induced from the deRham vector field $d_{dR}$, which maps $m$-forms to $(m+1)$-forms. If $Q$ is a polynomial curved $A_\infty$-structure on $V$  then $\Omega^\bullet_{cycl}(V) $ carries another differential induced from $Lie_Q$,  which preserves each space of $m$-forms.

We also need to consider the space of germs of differential forms at the origin in $V$. We define the space of germs of de Rham forms $\hat{\Omega}^\bullet(V)$ to be the completion of $\Omega^\bullet(V)$ at the ideal generated by $V$, so for example
$$\hat{\Omega}^1(V) = \hat{\cO}(V)\otimes V[1]^\vee \otimes \hat{\cO}(V)$$
The space of germs of cyclic forms $\hat{\Omega}_{cycl}^\bullet(V)$ is just the coinvariants of $\hat{\Omega}^\bullet(V)$ under the cyclic group action.

If $Q$ is a non-polynomial $A_\infty$-structure on $V$ then it gives only the germ of a vector field at the origin. Hence $Lie_Q$ is not defined on differential forms over the whole of $V$, but it is defined on germs of differential forms at the origin.

\subsubsection{The infinite-dimensional case}

Now we recall how to handle the situation when $V$ is not degree-wise finite-dimensional. In this case we want to avoid dualizing $V$, so instead of considering the free algebra $\cO(V)$ we instead consider the cofree coalgebra
$$\mathcal{C}(V):= \prod_{k\geq 0} V^{\otimes k}$$
If $V$ were degree-wise finite-dimensional, this would be the graded linear dual of $\cO(V)$. The coproduct is the \quotes{shuffle} coproduct, which takes a monomial to the sum of all ways of splitting it in two (e.g. \cite{keller1}). Notice that the cofree coalgebra is a direct product, not a direct sum. All of the above constructions go through in their dual version. For example we have a sequence of subcoalgebras
$$\bigoplus_{i=0}^k V^{\otimes i} \subset \mathcal{C}(V)$$
corresponding to $k$th order neighbourhoods of the origin, and if we take the colimit over this sequence we get the cocompleted coalgebra
$$\hat{\mathcal{C}}(V) := \bigoplus _{k\geq 0} V^{\otimes k}$$
which corresponds to a formal neighbourhood of the origin.

In this setting, a \textit{vector field} on $V$ is a coderivation of $\mathcal{C}(V)$, and the germ of a vector field at the origin is a coderivation of $\hat{\mathcal{C}}(V)$.

\begin{defn} Let $V$ be a graded vector space. A curved $A_\infty$-structure on $V$ is a coderivation
$$Q :\hat{\mathcal{C}}(V) \to \hat{\mathcal{C}}(V)$$
of degree 1, such that
$$[Q,Q] =0$$
\end{defn}
As before, $Q$ is determined by its Taylor coefficients, which are maps
$$Q_k : V^{\otimes k} \to V$$
and we say that the structure is \textit{polynomial} if only finitely many of these are non-zero. The induced structure on the degree-shifted vector space $A=V[-1]$ is called an $A_\infty$-algebra.
If $V$ is degree-wise finite-dimensional these definitions are equivalent to our previous ones, just by dualizing the Taylor coefficients.

The \textit{dual de Rham differential forms} on $V$ are given by the coalgebra
$$\mho^\bullet(V) := \mathcal{C}(V\oplus V[1])$$
It carries a coderivation $d'_{dR}$ defined by dualizing the definition of $d_{dR}$. The dual notion to cyclic functions is given by the cyclically invariant subspace
$$\mathcal{C}_{cycl}(V) = (\mathcal{C}(V))^{\Z}$$
and hence the \textit{dual cyclic differential forms} are given by the cyclically invariant dual de Rham forms
$$\mho^\bullet_{cycl}(V) := \mathcal{C}_{cycl}(V\oplus V[1])$$
This carries an induced differential $d'_{cycl}$ that maps $\mho^m_{cycl}(V)$ to $\mho^{m-1}_{cycl}(V)$.

The space of germs of dual de Rham forms at the origin is given by the co-completion of $\mho^\bullet(V)$ at the co-ideal generated by $V$, we denote it by $\hat{\mho}^\bullet(V)$. Taking cyclic invariants we get the space of germs of dual cyclic forms $\hat{\mho}_{cycl}^\bullet(V) = (\hat{\mho}^\bullet(V))^{\Z}$.

If $Q$ is a polynomial $A_\infty$-structure on $V$ it induces a coderivation $Lie_{Q}$ of $\mho^\bullet(V)$ and a differential on each $\mho^m_{cycl}(V)$, again just by dualizing the definitions from the previous section. If $Q$ is not polynomial, then $Lie_{Q}$ is only defined on the germs of dual differential forms at the origin.

\subsubsection{The Hochschild and Borel-Moore Hochschild complexes}

Let $(V,Q)$ be a polynomial curved $A_\infty$-structure, and $A=V[-1]$ the associated  $A_\infty$-algebra with products $m_k$.

\begin{defn} \cite{CT} The Borel-Moore Hochschild chain complex of $A$ is
$$C_\bullet^{\Pi}(A) := (\mho^1_{cycl}(V)[-2], Lie_{Q})$$
\end{defn}

If $Q$ is not polynomial, so $A$ has infinitely many non-zero products, then this is not defined. In that case we can instead take the germs of dual cylic 1-forms at the origin, and get the more classical  (e.g. \cite{lowen}) Hochschild chain complex:
\begin{defn} The Hochschild chain complex of $A$ is
$$C_\bullet(A) := (\hat{\mho}^1_{cycl}(V)[-2], Lie_{Q})$$
\end{defn}
 Recall that the space of de Rham dual 1-forms is given by
$$\mho^1(V) = \prod_{k\geq 0} V^{\otimes k} \otimes V[1]\otimes \prod_{l\geq 0} V^{\otimes l} $$
It follows that the graded vector space underlying the Borel-Moore Hochschild complex is
$$ \mho_{cycl}^1(V) [-2] \cong  V[-1]\otimes \prod_{k\geq 0} V^{\otimes k}  = A\otimes \prod_{k\geq 0} A^{\otimes k}[k]$$
 The space underlying the Hochschild complex is exactly the same, except that the direct product becomes a direct sum. The differential is the same on both complexes, to describe it explicitly we must unpack the definition of $Lie_Q$ and insert a lot of signs. We give the first three terms of the differential, this is sufficient to cover the case when $A$ is a curved dga. They are:
\begin{eqnarray*}d_2( a_0\otimes...\otimes a_{k})& =& \sum_{i=0}^{k-1} (-1)^{|a_0|+...+|a_i|+i+1}\, a_0\otimes...\otimes a_i
a_{i+1}\otimes...\otimes a_{k}\\
&& \;\;+ \;\;(-1)^{(|a_0|+...+|a_{k-1}| + k +1)(|a_k| + 1)}\, a_k a_0\otimes ...\otimes a_{k-1} \end{eqnarray*}
$$d_1( a_0\otimes...\otimes a_k) = \sum_{i=0}^{k} (-1)^{|a_0|+...+|a_{i-1}|+i}\, a_0\otimes...\otimes da_{i}\otimes...\otimes a_{k} $$
\begin{equation}\label{hochschilddiff}d_0(a_0\otimes..\otimes a_{k}) = \sum_{i=0}^{k} (-1)^{|a_0|+...+|a_{i}|+i+1}\, a_0\otimes...\otimes a_{i} \otimes W \otimes a_{i+1}\otimes...\otimes a_{k}\end{equation}

The complicated sign for the last term of $d_2$ arises from permuting $a_k$ through the other elements.

\subsubsection{Curved $A_\infty$-categories}

We can also use this geometric language to describe curved $A_\infty$-categories. Fix a set $Ob$ of objects, and pick a graded vector space $\mathcal{V}(a,b)$ for each ordered pair of objects $(a,b)$. This is the same thing as a graded bimodule $\mathcal{V}$ over the semi-simple algebra $\C^{Ob}$ generated by the objects. It's not very misleading to think of $\mathcal{V}$ as just a single vector space, so we can try and perform all the above constructions of non-commutative geometry on it. This works fine, as long as we remember that the ground ring is $\C^{Ob}$. For example, $\mathcal{V}\otimes \mathcal{V}$ must be read as a tensor product over $\C^{Ob}$, so $\mathcal{V}^{\otimes k}$ is the $\C^{Ob}$-bimodule with components
$$\mathcal{V}^{\otimes k}(a,b) = \bigoplus_{c_1,...,c_{k-1}\in Ob} \mathcal{V}(a, c_1)\otimes \mathcal{V}(c_1,c_2)\otimes... \otimes \mathcal{V}(c_{k-1}, b) $$
Now we can define vector fields, differential forms, etc. as before, and declare that a curved $A_\infty$-category with objects $Ob$ is a curved $A_\infty$-structure on a $\C^{Ob}$-bimodule $\mathcal{V}$.

Note that the cyclic invariants in $\mathcal{V}^{\otimes k}$ are not a $\C^{Ob}$-bimodule, they are the vector space
$$(\mathcal{V}^{\otimes k})^{\Z} = \left(\bigoplus_{a\in Ob} \mathcal{V}^{\otimes k}(a,a)\right)^\Z $$
So for example if $\mathcal{A}$ is a curved $A_\infty$-category, then its Borel-Moore Hochschild homology is a chain-complex with underlying graded vector space
$$\prod_{k\geq 0} \;\;\bigoplus_{a_0,...,a_k\in Ob(\mathcal{A})} \mathcal{A}(a_0,a_1)\otimes ...\otimes \mathcal{A}(a_k,a_0) [k]$$

Also, since a point of a vector space $V$ is the same as a linear map $\C\to V$, choosing a \quotes{point} of $\mathcal{V}$ means giving a map of $\C^{Ob}$-bimodules
$$\C^{Ob} \to \mathcal{V} $$
which is the same thing as choosing an element of $\mathcal{V}(a,a)$ for all $a\in Ob$.

\subsubsection{Translation maps}

Let $V$ be a graded vector space, and let $x\in V$ be an element. In ordinary commutative ungraded geometry $x$ would induce a constant vector field on $V$ and an affine endormorphism of 'translate by $x$'. In our non-commutative context, the constant vector field field is the coderivation
$$X: \mathcal{C}(V) \to \mathcal{C}(V)$$
with only a constant Taylor coefficient $x:\mathbb{C}\to V$. Hence we define
$$T_x := exp(X) = \sum_{n\geq 0}\frac{1}{n!}X^n: \mathcal{C}(V) \to \mathcal{C}(V) $$
 as the appropriate analogue of translation by $x$. This is a map of coalgebras, its components are
$$ \sum_{s_0+...+s_l = k} \id^{\otimes s_0}\otimes x\otimes \id^{\otimes s_1}\otimes...\otimes x\otimes \id^{\otimes s_l} \;:\;\; V^{\otimes k} \to V^{\otimes k+l} $$
More generally, if $\mathcal{V}$ is a $\C^{Ob}$-bimodule over a set of objects $Ob$, then we have a translation map $T_x$ on $\mathcal{V}$ for each element $x: \C^{Ob} \to \mathcal{V}$.

\subsection{Borel-Moore Hochschild complexes of categories of cdg-modules}

We now return to our main object of study: the category of perfect cdg-modules over a curved dga. We want to understand the Borel-Moore Hochschild homology of this category, and we will approach this using the geometric language that we have been setting up in the previous sections.

 Let $(A,W)$ be a curved dga, and let
$$\mathcal{P} \subset \Perf(A,W)$$
be a full sub-category of the category of perfect $(A,W)$-cdg-modules, so $\mathcal{P}$ is a dg-category. Let $\tilde{\mathcal{P}}$ be the underlying graded category of $\mathcal{P}$ (i.e. forget the differential), so $\tilde{\mathcal{P}}$ is a full subcategory
$$\tilde{\mathcal{P}} \subset proj(A)$$
of the category of finitely-generated projective graded $A$-modules.

 Let $\mathcal{V}$ be the $Ob(\mathcal{P})$-bimodule underlying $\mathcal{P}[1]$. The dg-category structure on $\mathcal{P}$ is encoded in a vector field
$$Q = Q_2 + Q_1: \mathcal{C}(\mathcal{V}) \to \mathcal{C}(\mathcal{V}) $$
which has only quadratic and linear Taylor coefficients. If we use only the quadratic term $Q_2$ then we are encoding the graded category $\tilde{P}$.

The linear term $Q_1$ is of a particular form. Consider the following (degree zero) element of $\mathcal{V}$
$$d: \id_M \mapsto d_M $$
There is an associated constant vector field $D$ on $\mathcal{V}$. Then
$$Q_1 = [Q_2, D] $$
This is just by the definition (\ref{dhom}) of the differential in $\Perf(A,W)$.  We can also define a constant vector field
$$Q_0 := \frac12 [[Q_2.D],D]$$
which corresponds to the element
$$Q_2(D\otimes D) : \id_M \mapsto  W\id_M $$
This is the curvature term from Remark \ref{curvedcat}. Each of these three terms $Q_2, Q_1, Q_0$ are degree 1, and they commute with each other and themselves, so any combination of them  encodes a curved $A_\infty$-structure. Two of these possible structures are given by $\mathcal{P}$ and $\tilde{\mathcal{P}}$, the third one that we need to consider is the one encoded by the vector field $Q_2 - Q_0$. This is a curved category, it's obtained from the graded category $\tilde{\mathcal{P}}$ by adding in a curvature term given by the central element
$$ \id_M \mapsto -W\id_M$$
We'll denote this curved category by $(\tilde{\mathcal{P}}, -W)$.

Now consider the \quotes{translation by $d$} map
$$T_d: \mathcal{C}(\mathcal{V}) \to \mathcal{C}(\mathcal{V})  $$
Since it is an isomorphism (the inverse is $T_{-d}$) we can ask what effect it has
on a vector field $Y$. It is elementary that
\begin{equation*}
T_{-d} Y T_{d} = Y + [Y,D] + \frac12 [[Y,D],D] + \frac16 [[[Y,D],D],D] + ...
\end{equation*}
One way to view this formula is as giving the relationship between the Taylor expansions of
Y at zero and at the point $d\in\mathcal{V}$. Let's apply this formula to the case $Y=Q_2$. Since this is a quadratic vector field, only the first three terms on the RHS are non-zero, and we have
$$T_{-d} Q_2 T_{d} = Q_2 + Q_1 + Q_0 $$
Also $Q_0$ is a constant vector field, so $T_{-d} Q_0 T_{d} = Q_0$, and hence
$$T_{-d}(Q_2 - Q_0) T_d = Q_2 + Q_1 $$
So the translation map $T_{d}$ gives us a polynomial isomorphism between $\mathcal{P}$ and  $(\tilde{\mathcal{P}}, -W)$.

\begin{lem}\label{Tdiso} We have an isomorphism
 $$T_d: C^{\Pi}_\bullet(\mathcal{P}) \iso C^{\Pi}_\bullet(\tilde{\mathcal{P}}, -W) $$
between the Borel-Moore Hochschild complexes of $\mathcal{P}$ and $(\tilde{\mathcal{P}}, -W)$.
\end{lem}
\begin{proof}
 $T_d$ is an isomorphism of the non-commutative affine space $\mathcal{V}$ that intertwines the two vector fields $Q_2 + Q_1$ and $Q_2 - Q_0$. Therefore it induces an isomorphism
$$(\mho^1_{cycl}(\mathcal{V}), Lie_{Q_2 + Q_1}) \iso (\mho^1_{cycl}(\mathcal{V}), Lie_{Q_2 - Q_0})$$
\end{proof}

To write $T_d$ explicitly, let $(M_0, d_0),...,(M_k,d_k)$ be cdg-modules in $\mathcal{P}$, and let
$$M_0 \stackrel{\alpha_0}{\to} M_1 \stackrel{\alpha_1}{\to}... \stackrel{\alpha_{k-1}}{\to} M_k \stackrel{\alpha_k}{\to} M_0$$
be morphisms. Then
$$T_d:\alpha_0\otimes...\otimes\alpha_k \mapsto \sum_{s_0,...,s_k\geq 0}  \alpha_0 \otimes (d_1)^{\otimes s_1} \otimes \alpha_1 \otimes (d_2)^{\otimes s_2}\otimes... \otimes \alpha_k \otimes (d_0)^{\otimes s_0}  $$
Note that $\mathcal{P}$ and $(\tilde{\mathcal{P}}, -W)$ are not $A_\infty$-quasi-isomorphic in the conventional sense, because this map $T_d$ between them doesn't preserve the origin in $\mathcal{V}$. Consequently we don't get a map on the ordinary Hochschild complexes, because the ordinary Hochschild complex is the cocompletion of the Borel-Moore Hochschild complex at the origin. Instead, $T_d$ maps $C_\bullet(\mathcal{P})$ to a different cocompletion of $C^\Pi_\bullet(\tilde{\mathcal{P}}, -W)$, not the cocompletion at the origin, but the cocompletion at the point $-d\in \mathcal{V}$.

It follows that Borel-Moore Hochschild homology is rather weaker than ordinary Hochschild homology (at least in this context), as it is always equal to the homology of $C^\Pi_\bullet(\tilde{\mathcal{P}}, -W)$ so it doesn't depend on the differential in $\mathcal{P}$.  For example, take $A$ to be an ordinary algebra, and let $\mathcal{P}$ be the subcategory of $\Perf(A)$ containing the single dg-module
$$A \stackrel{a}{\to} A[-1]$$
for some element $a\in A$. Then the Borel-Moore Hochschild homology is independent of $a$, but the Hochschild homology ranges from zero (when $a$ is a unit) to $HH_\bullet(A)$ (when $a=0$).

\subsubsection{The generalized trace map}\label{tracemap}

Let $A$ be a graded algebra (with no differential), and let
$$\tilde{\mathcal{P}} \subset proj(A)$$
be a full subcategory. Our first result in this section is that the Hochschild complexes of $\mathcal{\tilde{P}}$ and $A$ are quasi-isomorphic. This is well-known, and we only include it because we want to know the quasi-isomorphism explicitly, and we could not find the formulas in the literature for the graded case. Our strategy is to embed $\tilde{\mathcal{P}}$ in an (infinite-rank) matrix algebra over $A$, and then use the proof from \cite{loday} with some additional signs.

We need to assume:
\begin{equation}\label{containsA}\mbox{ There is a module $N\in\tilde{\mathcal{P}}$ that contains $A$ as a direct summand.} \end{equation}
Given such an $N$ we have a linear embedding
$$\iota: A \into \Hom_A(N,N) $$
which induces a chain map on the Hochschild complexes
$$\iota: C_\bullet(A) \into C_\bullet(\Hom_A(N,N)) \into C_\bullet(\tilde{\mathcal{P}})$$
We want to write down a homotopy inverse to $\iota$. To do this we first fix, for every $M\in \tilde{\mathcal{P}}$, an embedding
$$M \into \bigoplus_{i=1}^{r_M}A[\sigma_i]$$
 of $M$ as a direct summand of a graded, finite-rank free $A$-module. This means that any morphism
$$\alpha\in \Hom_{\tilde{\mathcal{P}}}(M,M')$$
is explicitly a matrix over $A$ of size $r_M\times r_{M'}$. As usual let's write
$$\alpha_{ij}: A[\sigma_i] \to A[\sigma_j]$$
for the entries of the matrix (our matrices act on the right). If $\alpha$ is homogeneous we have
$$|\alpha| = |\alpha_{ij}| + \sigma_i - \sigma_j $$
for each $i,j$. Let's also write
$$\alpha_{ i \bullet } : A[\sigma_i] \to M' $$
$$\alpha_{ \bullet j} : M\to A[\sigma_j] $$
for the maps given by the rows and columns of the matrix. Finally, for any $M$ we denote by $\epsilon_i$ the composition
$$\epsilon_i: N\to  A \to A[\sigma_i] \to M $$
This has degree $\sigma_i$.

\begin{defn} \label{Trdef}Let $M_0,...,M_t\in\tilde{\mathcal{P}}$. We define generalized trace maps
$$\Tr : \Hom_A(M_0, M_1)\otimes...\otimes  \Hom_A(M_{t-1},M_t)\otimes \Hom_A(M_t, M_0) \to A^{\otimes (t+1)} $$
$$ \Tr(\alpha^0\otimes...\otimes \alpha^{t-1}\otimes \alpha^t) = \sum (-1)^\sigma \alpha^0_{i_0 i_1}\otimes...\otimes \alpha^{t-1}_{i_{t-1} i_t}\otimes \alpha^t_{i_t i_0} $$
 where every $i_k$ ranges from  $1$ to $r_{M_k}$. The sign is given by
 $$\sigma = (|\alpha^0_{i_0 i_1}| + ... + |\alpha^t_{i_t i_0}| + t+ 1)\sigma_{i_0} + \sigma_{i_1} +...+\sigma_{i_t}$$
\end{defn}
If $A$ is concentrated in even degrees, we can instead write the sign as
\begin{equation}\label{Trsign}\sigma = \sigma_{i_0} + |\alpha^1| + |\alpha^3|+ ... + |\alpha^{s}|
\end{equation}
where $s=t$ or $t-1$ according to whether $t$ is odd or even. In particular for $t=0$ the map is
$$\alpha \mapsto \sum_{i} (-1)^{\sigma_i } \alpha_{ii}$$
which is the standard supertrace map. 

Adding these maps together we get a linear map
$$\Tr: C_\bullet(\tilde{\mathcal{P}}) \to C_\bullet(A)$$

\begin{lem}\label{Triso} The map
$$Tr: C_\bullet(\tilde{\mathcal{P}}) \to C_\bullet(A)$$
is a chain map, and is homotopy inverse to $\iota$, so it is a quasi-isomorphism.
\end{lem}
\begin{proof}
We lift the proof from \cite[Thm. 1.2.4]{loday}, the only addition is the signs. Firstly, the differentials on each side are given by the expression for $d_2$  in the formulas (\ref{hochschilddiff}), and the map $\Tr$ commutes with every term in this expression, so it is a chain map. The composition $\Tr\circ \iota$ is the identity on $C_\bullet(A)$, so to complete the proof we just need a homotopy between $\iota \circ \Tr$ and the identity on $C_\bullet(\tilde{\mathcal{P}})$.

For
$$\alpha^0\otimes...\otimes \alpha^t \in \Hom_A(M_0, M_1)\otimes...\otimes\Hom_A(M_t, M_0)$$
we define, for each $s\in [0,t]$,
$$h_t^s(\alpha^0\otimes...\otimes \alpha^t)=\sum (-1)^\sigma \iota \alpha^0_{\bullet i_1}\otimes \iota(\alpha^1_{i_1 i_2})\otimes...\otimes \iota (\alpha^s_{i_s i_{s+1}})\otimes \epsilon_{i_{s+1}} \otimes \alpha^{s+1}\otimes...\otimes \alpha^t $$
This is an element of
\begin{multline*}
\Hom_A(M_0, N)\otimes \Hom_A(N,N)^{\otimes s} \otimes \Hom(N, M_{s+1})\otimes \Hom_A(M_{s+1}, M_{s+2}) \otimes ...\\ ...\otimes  \Hom_A(M_t,M_0)\end{multline*}
The sign is
$$\sigma = |\alpha^0| + ... + |\alpha^s| + \sigma_{i_1} + ... + \sigma_{i_{s+1}} $$
Now let
$$h_t = \sum_{s=0}^t (-1)^s h_t^s $$
This is a degree zero map, so adding all the $h_t$ together gives an endomorphism of $C_\bullet(\tilde{\mathcal{P}})$ of degree $-1$. This is the required homotopy.
\end{proof}

If we now pick a central degree two element $W\in A$, then the pair $(A,W)$ forms a curved algebra. The Hochschild complex of $(A,W)$ is obtained from the Hochschild complex of $A$ by adding in a new term to the differential, the $d_0$ term from (\ref{hochschilddiff}).

We can also use $W$ to turn $\tilde{\mathcal{P}}$ into a curved category, as we did in the previous section. The curvature is the central element
$$W: \id_M \mapsto W\id_M $$
The Hochschild complex of $(\tilde{\mathcal{P}}, W)$ is similarly a deformation of the Hochschild complex of $\tilde{\mathcal{P}}$ by a $d_0$ term. We still have a chain map
$$Tr: C_\bullet(\tilde{\mathcal{P}},W) \to C_\bullet(A,W)$$
because the map $Tr$ commutes with every summand of $d_0$. This is not a quasi-isomorphism, however we do get a quasi-isomorphism between the Borel-Moore Hochschild complexes, under an additional assumption.

Because $A$ is just a graded algebra, with no differential or curvature, the Hochschild complex of $A$ is actually bi-graded. The first grading is the internal grading on $A$, and the second is by the number of tensor powers of $A$. We need to assume that the Hochschild homology of $A$ is bounded with respect to this second grading, i.e. we are assuming that the ungraded algebra underlying $A$ has bounded Hochschild homology.

\begin{prop} \label{Trqi} Assume that the Hochschild homology of $A$ is bounded in the above sense. Then for any degree 2 central element $W\in A$, we have a quasi-isomorphism
$$Tr: C^\Pi_\bullet (\tilde{\mathcal{P}}, W) \to C^\Pi_\bullet (A,W)$$
\end{prop}
\begin{proof}
This argument is based very closely on \cite[Sect. 4.9]{CT}, so we will be brief and refer the reader to there or to \cite{loday} for a clearer explanation. Recall that the differential on $C^\Pi_\bullet (\tilde{\mathcal{P}}, W)$ is a sum of two commuting differentials $d_2$ and $d_0$. Also, note that we can split $C^\Pi_\bullet (\tilde{\mathcal{P}}, W)$ into a direct sum of the following two pieces:
$$C^\Pi_\bullet (\tilde{\mathcal{P}}, W)^{ev} := \prod_{k \makebox{\small{ even}}}\tilde{\mathcal{P}}^{\otimes k+1} [k] $$
$$C^\Pi_\bullet (\tilde{\mathcal{P}}, W)^{od} := \prod_{k \makebox{\small{ odd}}}\tilde{\mathcal{P}}^{\otimes k+1} [k] $$
The differential exchanges these two pieces, so we can consider $C^\Pi_\bullet (\tilde{\mathcal{P}}, W)$ to be $\mathbb{Z}\times \mathbb{Z}_2$ -graded.

Now let  $\mathbf{x}$ and $\mathbf{y}$ be formal variables, and consider the bi-graded vector space
$$BC(\tilde{\mathcal{P}}, W) = \bigoplus_{i,j \in \Z} \left( \tilde{\mathcal{P}}^{\otimes j-i +1}[j-i] \right) \mathbf{x}^{i} \mathbf{y}^{j}$$
The first grading is the internal grading on each $\tilde{\mathcal{P}}^{\otimes j-i +1}[j-i]$, and the second grading comes from giving both $\mathbf{x}$ and $\mathbf{y}$ bi-degree $(0,1)$. Shifting the second degree by 2 is an isomorphism. We equip $BC(\tilde{\mathcal{P}}, W) $ with the differential
$$\mathbf{y}^{-1} d_2 + \mathbf{x}^{-1} d_0 $$
We also consider the bi-complexes
$$Q_p BC(\tilde{\mathcal{P}}, W) = \bigoplus_{i\geq p}\bigoplus_{j \in \Z} \left(\tilde{\mathcal{P}}^{\otimes j-i +1} [ i-j]\right) \mathbf{x}^{i} \mathbf{y}^{j}$$
with the same differential. These are quotients of $BC(\tilde{\mathcal{P}}, W)$ by the subcomplex where $i<p$, and form a sequence :
\begin{equation} ... \to Q_p BC(\tilde{\mathcal{P}}, W)  \to Q_{p+1}BC(\tilde{\mathcal{P}}, W) \to... \end{equation}
Let $Q_{-\infty}BC(\tilde{\mathcal{P}}, W)$ be the (inverse) limit over this sequence. To see what this is, notice that the piece of $Q_p BC(\tilde{\mathcal{P}}, W)$ having second degree equal to $t$ is
\begin{eqnarray*}Q_pBC(\tilde{\mathcal{P}}, W)_{\bullet, t} &=& \bigoplus_{i \geq p} \left( \tilde{\mathcal{P}}^{\otimes t-2i +1}[t-2i] \right) \mathbf{x}^{i} \mathbf{y}^{t-i}\\
&\cong & \bigoplus_{\substack{k \leq t-2p \\ k\equiv t (mod 2)}}\tilde{\mathcal{P}}^{\otimes k +1}[k]
\end{eqnarray*}
Taking the limit $p\to -\infty$ we get that $Q_{-\infty}BC(\tilde{\mathcal{P}}, W)_{\bullet, t}$ is either $C^\Pi_\bullet (\tilde{\mathcal{P}}, W)^{ev}$ or $C^\Pi_\bullet (\tilde{\mathcal{P}}, W)^{od}$ depending on whether $t$ is even or odd. If we quotient by the 2-periodicity in $t$ we get back $C^\Pi_\bullet (\tilde{\mathcal{P}}, W)$ with its $\Z\times \Z_2$-grading.

We can perform exactly the same constructions with $(A,W)$, and we have in particular chain maps
\begin{equation}\label{Qp}Tr: Q_pBC(\tilde{\mathcal{P}}, W) \to Q_pBC(A, W)\end{equation}
If we draw either of these bi-complexes in the $\mathbf{x}$-$\mathbf{y}$ plane we see that we can compute their homology using a spectral sequence. If we take the homology of $d_2$ first, then on page 1 we get an infinite number of copies of the Hochschild homology of either $\tilde{\mathcal{P}}$ or $A$. By Lemma \ref{Triso}, the map $Tr$ induces an isomorphism between page 1 on either side. Furthermore, by our boundedness assumption all terms on page 1 are zero if we move far enough up from the diagonal, so the spectral sequences eventually collapse, and we deduce that (\ref{Qp}) is a quasi-isomorphism.

Our boundedness assumption also implies that when the second degree $t$ is large enough the homology of $Q_pBC(A, W)$ is independent of $p$. This, together with the surjectivity of
$$Q_pBC(A, W) \to Q_{p+1}BC(A, W)$$
is enough to guarantee that the homology of $Q_{-\infty}BC(A,W)$ is the limit of the homologies of the $Q_pBC(A, W)$. The same is true for $(\tilde{\mathcal{P}}, W)$, so we deduce that
\begin{equation*} Tr: Q_{-\infty}BC(\tilde{\mathcal{P}}, W) \to Q_{-\infty}BC(A, W)\end{equation*}
is a quasi-isomorphism. Quotienting by the 2-periodicity in $t$ we get the statement of the Proposition.
\end{proof}

Our boundedness condition on the Hochschild homolgy of $A$ is a kind of smoothness condition. Let's assume it holds, and fix a curvature element $W\in A$. Let
$$\mathcal{P}\subset \Perf(A,W)$$
be a full dg-subcategory, and let $\tilde{P}\subset proj(A)$ be the underlying graded category. We assume that $\tilde{P}$ satisfies (\ref{containsA}), but this is a very weak assumption, as we can always achieve it by adding to $\mathcal{P}$ the contractible cdg-module
$$N = A\oplus A[-1] \;\;\;\;\;\;\;\;\;\; d_N = \left(\begin{array}{cc} 0 & W \\ 1 & 0 \end{array}\right) $$
In particular it holds when $\mathcal{P}$ is the whole of $\Perf(A,W)$. Then combining Lemma \ref{Tdiso} and Prop. \ref{Trqi}, we get:

\begin{thm}\label{mainthm} We have a quasi-isomorphism
$$Tr\circ T_d: C^\Pi_\bullet(\mathcal{P}) \iso C^\Pi_\bullet (A,-W) $$
\end{thm}

Using different methods, \cite{PolPos} have independently proved that these two complexes are quasi-isomorphic.

\section{Landau-Ginzburg B-models}

\subsection{The closed state-space}

\begin{defn} An affine Landau-Ginzburg B-model is the following data:
\begin{itemize}
\item A smooth $n$-dimensional affine variety $X$ over $\C$.
\item A choice of function $W\in \mathcal{O}_X$ (the \quotes{superpotential}).
\item An action of $\C^*$ on $X$ (the  \quotes{vector R-charge}).
\end{itemize}
such that
\begin{enumerate}
\item $-1\in \C^*$ acts trivially.
\item $W$ has weight (\quotes{R-charge}) equal to 2.
\end{enumerate}
\end{defn}

This means that $\mathcal{O}_X$ is a regular commutative algebra graded by the even integers, and $W\in \mathcal{O}_X$ is an
element of degree 2. Thus $(\mathcal{O}_X,W)$ is a curved algebra (with no odd graded part).

There is a weaker definition of vector R-charge where we keep only the trivial action of the
 the subgroup $\Z_2\subset \C^*$. This corresponds to working with $\Z_2$-graded
 curved algebras.

\begin{defn} The (off-shell) \textit{closed state space} of an affine LG B-model $(X,W)$ is the graded vector space
$$\bigoplus \Omega_X^k [-k] $$
of holomorphic forms on $X$, with differential
$$ \alpha \mapsto dW\wedge \alpha $$
\end{defn}

Note that since $W$ has degree 2 the total degree of the differential is indeed 1.

If $W$ has an isolated singularities then the homology of this complex is
$$\Omega^{n}_X[-n]\; /\; (dW)$$
which is the Jacobi ring (times a volume form). It is well known in the
physics literature that this is the space of physical closed states.

\begin{defn} The category of B-branes for an affine LG model $(X,W)$ is the dg-category
$$Br(X,W) := \Perf(\mathcal{O}_X,W) $$
of perfect cdg-modules over the curved algebra $(R, W)$. Of course a finitely-generated projective $\mathcal{O}_X$-module is exactly a finite-rank vector
bundle on $X$.
\end{defn}

Since $\mathcal{O}_X$ has no odd graded part it follows that a brane $M$ splits as a direct sum $M_{ev}\oplus M_{od}$ where $M_{ev}$ (respectively $M_{od}$) is the sum of all factors that are shifted by an even (respectively odd) integer, and $d_M$ exchanges these two factors.

If $X=\C^n$, and we work with a $\Z_2$-grading, then since all vector bundles are trivial a B-brane $M$ is described by a pair of polynomial matrices $d_M^{od}$
and $d_M^{ev}$ such that  $d_M^{ev}d_M^{od}= d_M^{od}d_M^{ev} = W\id$. This
is a \quotes{matrix factorization} of $W$.

As discussed in the introduction, $Br(X,W)$ is supposed to be the open sector of the B-model TCFT constructed from $(X,W)$. If this is true, the canonical closed sector would be the Hochschild complex of $Br(X,W)$, and we expect that the closed state space $(\Omega_X, \wedge dW)$ is quasi-isomorphic to $C_\bullet(Br(X,W))$. However, all we can prove is that it is quasi-isomorphic to the Borel-Moore Hochschild complex of $Br(X,W)$, as we will now show.

In Theorem \ref{mainthm} we constructed a quasi-isomorphism between the Borel-Moore Hochschild complexes of $Br(X,W)$ and of the curved algebra $(\mathcal{O}_X, -W)$. To get to the closed state space, we use the map
\begin{equation}\label{phi}\phi:\mathcal{O}_X^{\otimes k+1} \to \Omega^k_X \end{equation}
$$ f_0\otimes...\otimes f_k \mapsto \frac{1}{k!} f_0 df_1\wedge...\wedge df_k $$
Using the definition (\ref{hochschilddiff}) of the Hochschild differential, it is elementary to check that $\phi$ gives a chain map from the (either standard or Borel-Moore) Hochschild complex of $(\mathcal{O}_X, -W)$ to the closed state space of $(X,W)$. When $W=0$, the famous theorem of Hochschild, Kostant and Rosenberg \cite{HKR} says that $\phi$ is a quasi-isomorphism
 $$\phi: C_\bullet(\mathcal{O}_X) \iso \Omega_X $$
When $W\neq 0$, C\u ald\u araru and Tu have shown \cite[Thm. 4.2]{CT} that it instead gives a quasi-isomorphism
$$\phi: C^\Pi_\bullet(\mathcal{O}_X, -W) \iso (\Omega_X , dW\wedge)$$
from the Borel-Moore Hochschild complex. Note that for their theorem they assume that $W$ has isolated singularities, however if one wants only this statement then that is unnecessary because the boundedness of $\Omega_X$ alone causes the degeneration of the relevant spectral sequence (see the proof of Prop. \ref{Trqi}).
\begin{cor} \label{maincor}The Borel-Moore Hochschild complex of $Br(X,W)$ is quasi-isomorphic to the closed state space of $(X,W)$, under the map
$$\phi\circ \Tr\circ T_d: C_\bullet^\Pi(Br(X,W)) \iso \Omega_X^\bullet $$
\end{cor}
We can write this quasi-isomorphism explicitly by unpacking Thm. \ref{mainthm}. Let $M_0,...,M_k$ be B-branes, with curved differentials $D_0,...,D_k$, and suppose we have morphisms
$$M_0 \stackrel{\alpha_0}{\to} M_1 \stackrel{\alpha_1}{\to}... \stackrel{\alpha_{k-1}}{\to} M_k \stackrel{\alpha_k}{\to} M_0$$
each of homogeneous degree. In order that the generalized trace map is defined, we have to explicitly write each $M_i$ as a summand of a trivialized free vector bundle, so each $D_i$ and $\alpha_i$ is a matrix of elements of $\mathcal{O}_X$. Then the element 
$$\alpha_0\otimes...\otimes \alpha_k \in C^\Pi_\bullet(Br(X,W))$$
 maps to 
$$\sum_{s_0,..., s_k \geq 0} \Tr \left( \alpha_0 \otimes (D_1)^{\otimes s_1} \otimes \alpha_1 \otimes (D_2)^{\otimes s_2}\otimes... \otimes \alpha_k \otimes (D_0)^{\otimes s_0}  \right) $$
in $C^\Pi_\bullet(\mathcal{O}_X, -W)$. Applying $\phi$ to this gives us
\begin{equation}\label{higherbb}\sum_{s_0,..., s_k \geq 0} \frac{(-1)^{\tau}}{(k+s_0+...+s_k)!}\; \mbox{Tr} \Big( \;\alpha_0 (dD_1)^{ s_1}  d\alpha_1  (dD_2)^{ s_2}...  d\alpha_k  (dD_0)^{ s_0}  \Big) \end{equation}
in $\Omega_X$. When we write $dD_i$ or $d\alpha_i$ here we mean the matrices of one forms obtained by applying $d$ to each entry in $D_i$ or $\alpha_i$, we are then multiplying these matrices together (over the ring $\Omega_X$) and taking the supertrace. To get this expression we are using the form (\ref{Trsign}) of the sign in Tr, which is valid since $\mathcal{O}_X$ is concentrated in even degrees. The sign $(-1)^\tau$ is given by
$$\tau = \left\lceil \frac{k+ \sum s_i}{2}\right\rceil + \sum_{i=1}^k (s_1+...+s_i + i)(|\alpha_i| +1)$$
The Hochschild complex is a sub-complex of the Borel-Moore Hochschild complex, so by restriction we have a chain map
\begin{equation}\label{chainonHH}\phi\circ \Tr\circ T_d: C_\bullet(Br(X,W)) \to \Omega_X\end{equation}

It should follow from \cite[Cor. B]{PolPos} that the inclusion of the Hochschild complex into the Borel-Moore Hochschild complex is a quasi-isomorphism, provided that one can verify the necessary generation condition on the category. This would imply that (\ref{chainonHH}) is also a quasi-isomorphism. In the $\Z_2$-graded case the generation condition holds \cite{LinPom}, but we must additionally require that the singular locus of $W$ is contained in $W^{-1}(0)$.

\subsection {The Kapustin-Li formula}
\label{KLformula}

Suppose that $W$ has isolated singularities. In this case $Br(X,W)$ should be a Calabi-Yau dg-category, and hence give the open sector of a TCFT, as we discussed in the introduction. This is rather delicate, it means that for all $M,N\in Br(X,W)$ we have a closed pairing
$$\Hom(M, N)\otimes \Hom(N,M) \to \C $$
which is symmetric and non-degenerate on homology. It has been known for a long time that the homotopy category of $Br(X,W)$ admits a non-degenerate pairing because of Auslander-Reiten duality \cite{auslander}, but this argument is non-constructive and doesn't give a chain-level pairing. 

Using path-integral methods, Kapustin and Li \cite{KL} derived a formula for a \quotes{trace map}
$$\End(M) \to \C$$ 
for any B-brane $M\in Br(X,W)$. This induces a chain-level pairing, and Dyckerhoff and Murfet have shown \cite{DycMur} that this is homologically non-degenerate. Unfortunately, it is not symmetric at the chain level.

\begin{figure}
 \subfigure[]{\label{disccorrelator}\includegraphics[width=.5\textwidth]{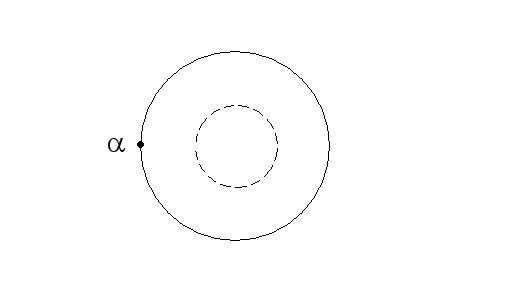}} \subfigure[]{\label{4pointdisc}\includegraphics[width=.5\textwidth]{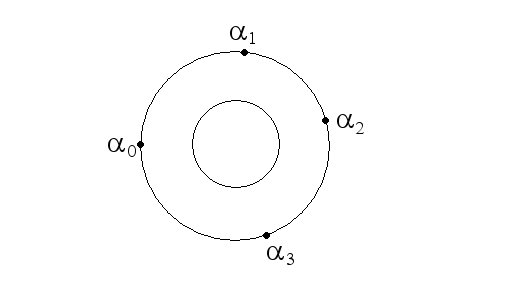}} 
\caption{Correlators over discs}
\label{relations}
\end{figure}

Physically, what Kapustin and Li compute is the 1-point correlator of an open string state $\alpha$ inserted on the boundary of a disc, as in Fig. \ref{disccorrelator}. We can factor this into two stages: if we cut the disc into an annulus and a smaller disc as indicated, then we can firstly propagate $\alpha$ to a closed string state living on the inner boundary of the annulus, and then take the correlator over the smaller disc.

The propagator over the annulus is called the boundary-bulk map. If assume that the closed state space is the Hochschild complex of $Br(X,W)$ then this map is tautological, it's the inclusion of $\End(M)$ into the Hochschild complex.

 More generally, there are \quotes{$n$-point boundary-bulk maps}, which are the propagators over annuli with $n$ open states inserted on the boundary, as in Fig. \ref{4pointdisc}. There is a slight subtley here: in these propagators we are varying the complex structure on the worldsheet, but only over a cell of dimension ${n-1}$ in the moduli space of complex structures \cite{costello2}.  Again, if the closed state space is the Hochschild complex then these maps are tautological, they map a composable set of morphisms $\alpha_0, ...., \alpha_{n-1}$ to the Hochschild chain $\alpha_0\otimes....\otimes \alpha_{n-1}$. To get something non-tautological, we can apply our chain map from Cor. \ref{maincor} to these Hochschild chains, and get elements in $\Omega_X$. 

The second stage is take the correlator over the disc with a closed state inserted on the boundary. Vafa \cite{vafa} argued that this is given by the residue map
$$\mbox{Res}_W:  \Omega_X \to \C $$
which takes a form to the sum of its residues at the singularities of $W$. This vanishes on the image of $\wedge dW$ (e.g. \cite[III.9]{hartshorne}), so it gives a closed element of the dual of the closed state space.

\begin{cor}We have a closed element of the dual of $C_\bullet(Br(X,W)$ (i.e. a Hochschild cocycle), given by $Res_W\circ \phi\circ \Tr\circ T_d$.\end{cor}

Applying this to a single endomorphism $\alpha$ of a B-brane $(M, D)$, we get (using (\ref{higherbb}))
\begin{equation}\label{formulafordisc}\langle\alpha\rangle_{disc} = \frac{(-1)^{\lceil\frac{n}{2}\rceil}}{n!} \mbox{Res}_{W}\left(  \mbox{Tr} \big( \;\alpha (dD)^ n  \big)\right) \end{equation}
 Kapustin and Li  work in the case $X=\C^n$, where the residue map can be written as a contour integral
$$ \mbox{Res}_W(\omega) = \frac{1}{(2 \pi i)^n} \oint \frac{\omega}{\prod_i \partial_i W} $$
where the integral is taken over a union of Lagrangian tori enclosing the singularities of $W$. In this case (\ref{formulafordisc}) becomes their formula, except that we have a correction to the sign. This agrees with the sign-correction found in \cite{DycMur}.

Thus we have recovered the Kapustin-Li trace map as the lowest-order term of a Hochschild cocycle on $Br(X,W)$. This solves, in some sense, the issue with chain-level symmetry, by use of the technology from \cite[Sect. 10]{KS} and \cite{ChoLee} (see also \cite[Sect 5.2]{DycMur} for a discussion of this point). Any negative cyclic cocycle (and in particular, any Hochschild cocycle) gives rise to a cyclic two-form on the non-commutative space underlying $Br(X,W)$. This two-form is called \quotes{symplectic} if its constant part gives a homologically non-degenerate pairing, which is true in our case by Dyckerhoff and Murfet's theorem. A symplectic two-form is the homotopy-invariant notion of a cyclic Calabi-Yau pairing, and furthermore any symplectic form can be made constant by an appropriate $A_\infty$-automorphism. We conclude that there is an alternative (but equivalent) $A_\infty$-structure on $Br(X,W)$ with respect to which the Kapustin-Li formula defines a cyclic Calabi-Yau pairing. In \cite{carqueville} the problem of explicitly determining this structure is addressed.

\subsection {Orbifolds}

We obtain more interesting and important examples of affine Landau-Ginzburg B-models if we allow the underlying space $X$ to be an orbifold, i.e. we
take a quotient
stack
$$X = [Y/G]$$
where $G$ is a finite group acting on a smooth affine variety $Y$, and add a superpotential
$W$ which is a $G$-invariant function on $Y$. In this section
we show how to adapt our results to this setting.

The natural definition of the category of B-branes on $(X,W)$ is the category of $G$-equivariant B-branes on $(Y,W)$, but we can recast this.
Recall that the \textit{twisted
group ring}
$$ A := \mathcal{O}_Y\rtimes \C[G] $$
is the vector space $\mathcal{O}_Y\otimes \C[G]$ with multiplication
$$(y_1\otimes g_1)\circ (y_2\otimes g_2) = (y_1g_1(y_2)\otimes g_1g_2)$$
$A$ inherits a grading and a superpotential $W\otimes 1$ from $R$, making it
a non-commutative curved algebra (the curvature is central since $W$ is invariant).
It is elementary to show that
 $$Br(X,W) := \Perf(A,W) $$
\begin{eg}
Let $X= [\C^2 / \Z_2]$ where $\Z_2$ acts with weight 1 on each co-ordinate.
We define a $\C^*$ R-charge action by letting $\C^*$ also act with weight 1 on each co-ordinate. Notice that $-1\in \C^*$ does indeed act trivially
on the orbifold (although not on $\C^2$). Let $x$ and $y$ be the two co-ordinates,
and let $W = x^2-y^2$.

The twisted group ring is $A = \C[x,y]\rtimes \C[\Z_2]$. Let $\tau$ be the
generator of $\Z_2$, then we have a complete pair of orthogonal idempotents
$$e_0 = \frac12 (1+\tau)\;\;\;\;\;\;\; e_1 = \frac12 (1-\tau)\;\;\;\;\;\;\;  e_i e_j = \delta_{ij}$$
This means we can write $A$ as a quiver algebra (with relations) using $e_1,
e_2$ as nodes. Every equivariant vector bundle on $\C^2$ is a direct sum
of the two line bundles $\mathcal{O}$ and $\mathcal{O}(1)$ associated to the
two characters of $\Z_2$. These correspond to the projective $A$-modules
$Ae_0$ and $Ae_1$. One example of a brane is given by $\mathcal{O}\oplus \mathcal{O}(1)$ with endomorphism
$$\left( \begin {array}{cc}
0&x+y\\
\noalign{\medskip}
x-y&0
\end {array}
\right)$$
This corresponds to the trivial $A$-module $A$ with endomorphism
$$e_0(x+y)e_1 + e_1(x-y)e_0 = x -y\tau$$
\end{eg}

Theorem \ref{mainthm} gives us a quasi-isomorphism between $C^\Pi_\bullet (Br(X,W))$
and $C^\Pi_\bullet(A, -W)$. As in the non-orbifold case, we can further map to a more geometric model for this complex.

For each $g\in G$, we can consider the fixed locus $Y^g$, and the restriction $W_g$ of $W$ to it. 

\begin{defn} The (off-shell) \textit{closed state space} of the affine orbifold
LG B-model $([Y / G], W)$ is the chain-complex of coinvariants
$$\left(\bigoplus_{g\in G} (\Omega^\bullet_{Y^g},  dW_g\wedge) \right)\; /
\;G$$
\end{defn}

There is a chain map $\psi$ from $C^\Pi_\bullet(A, -W)$ to this closed state space, defined as follows. We first map
$$(y_0\otimes g_0) \otimes (y_1\otimes g_1)\otimes....\otimes(y_k\otimes g_k) \in A^{\otimes k+1} $$
to 
$$y_0 \otimes g_0(y_1)\otimes... \otimes g_0g_1...g_{k-1}(y_k) \in \mathcal{O}_Y^{\otimes k+1
} $$
We map this to $\mathcal{O}_{Y^g}^{\otimes k+1}$ by restriction, and then to $\Omega_{Y^g}$ using the map $\phi$ from (\ref{phi}). We get $\psi$ by taking the direct sum over $g$ and taking coinvariants.

When $W=0$, Baranovsky \cite{baranovsky} has shown this is a quasi-isomorphism from the Hochschild complex $C_\bullet(A)$, and when $W\neq 0$ C\u ald\u araru and Tu \cite{CT}have shown that it's a quasi-isomorphism from the Borel-Moore Hochschild complex. 

\begin{cor} The Borel-Moore Hochschild complex of $Br([Y/G],W)$ is quasi-isomorphic to the closed state space of $([Y/G],W)$, under the map $\psi\circ \Tr\circ T_d$.
\end{cor}

If $W$ has isolated singularities, we can define a residue map on the closed state space by taking the sum over $g$ of 
$$\mbox{Res}_{W_g}: \Omega^\bullet_{Y^g} \to \C $$
since this is clearly well-defined on the coinvariants. Then as in the non-orbifold case we get a Hochschild cocycle on $Br([Y/G], W)$.

% ----------------------------------------------------------------

\bibliographystyle{amsplain}
\bibliography{mybib}

\end{document}